\documentstyle[11pt,dgspp]{article}

\newtheorem{theorem}{Theorem}[section]
\newtheorem{proposition}[theorem]{Proposition}

\newtheorem{corollary}[theorem]{Corollary}
\newremark{definition}[theorem]{Definition}
\newremark{example}[theorem]{Example}
\newremark{remark}[theorem]{Remark}

\newcommand{\bs}{\backslash}
\newcommand{\bul}{\,\scriptscriptstyle\bullet}
\newcommand{\dg}{\dot\gamma}

\newcommand{\ds}{\oplus} 
\newcommand{\g}{\gamma}

\newcommand{\half}{\mbox{$\txs\frac{1}{2}$}}
\newcommand{\io}{\iota}
\newcommand{\n}{\frak{n}}
\newcommand{\om}{\omega}
\newcommand{\sC}{\escr C}
\newcommand{\surj}{\rightarrow\kern-.82em\rightarrow}
\newcommand{\txs}{\textstyle}

\newcommand{\wO}{\widetilde{O}}
\newcommand{\z}{\frak{z}}
\newcommand{\E}{\frak{E}}
\newcommand{\G}{\Gamma}
\newcommand{\Iaut}{I^{a\kern-.05em u\kern-.02em t}}
\newcommand{\Ispl}{I^{s\kern-.055em p\kern-.025em l}}
\newcommand{\U}{\frak{U}}
\newcommand{\V}{\frak{V}}
\newcommand{\Z}{\frak{Z}}
\newcommand{\bZ}{\Bbb{Z}}
\newcommand{\A}{\Bbb A}
\newcommand{\B}{\Bbb B}
\newcommand{\R}{\Bbb{R}}

\renewcommand{\l}{\lambda}
\renewcommand{\v}{\frak{v}}

\makeatletter
\newcommand{\Ad}[1]{\mathop{\operator@font Ad}\nolimits_{#1}}
\newcommand{\ad}[1]{\mathop{\operator@font ad}\nolimits_{#1}}
\newcommand{\add}[2]{\mathop{\operator@font ad}\nolimits^{\dagger}_{#1}{#2}}
\newcommand{\Aut}{\mathop{\operator@font Aut}\nolimits}
\newcommand{\diag}{\mathop{\operator@font diag}\nolimits}   
\newcommand{\End}{\mathop{\operator@font End}\nolimits}
\newcommand{\mcp}{\mathop{\operator@font mult}\nolimits_{cp}}
\newcommand{\mev}{\mathop{\operator@font mult}\nolimits_{ev}}
\newcommand{\Ric}{\mathop{\operator@font Ric}\nolimits}
\newcommand{\Dspecp}{\mathop{\operator@font {\escr D}spec}\nolimits_\wp}
\newcommand{\specl}{\mathop{\operator@font spec}\nolimits_\ell}
\newcommand{\specp}{\mathop{\operator@font spec}\nolimits_\wp}
\newcommand{\tr}{\mathop{\operator@font tr}\nolimits}
\makeatother

\preprint{P9}
\title{Pseudoriemannian Nilpotent Lie Groups}
\author{Phillip E. Parker}
\address{Mathematics Department\\
     Wichita State University\\
     Wichita KS 67260-0033\\
     USA\\
     phil@math.wichita.edu}
\date{30 June 2005\\
      rev 3 October 2005 }

\abstract{This is a survey article with a limited list of references (as 
required by the publisher) which appears in the {\it Encyclopedia of
Mathematical Physics,} eds. J.-P. Fran\c{c}oise, G.\,L. Naber and Tsou
S.\,T.  Oxford: Elsevier, 2006. vol.\,4, pp.\,94--104.}

\msc{53C50}{22E25, 53B30, 53C30}

\begin{document}

\maketitle

\setcounter{page}{0}\thispagestyle{empty}\strut\vfill\eject

\section{Nilpotent Lie Groups}\label{nlg}
While there had not been much published on the geometry of nilpotent Lie
groups with a left-invariant Riemannian metric in 1990, the situation is
certainly better now; see the references in \cite{E}.  However, there is
still very little extant about the more general pseudoriemannian case.  In
particular, the 2-step nilpotent groups are nonabelian and as close as
possible to being Abelian, but display a rich variety of new and
interesting geometric phenomena \cite{CP}.  As in the Riemannian case,
one of many places where they arise naturally is as groups of isometries
acting on horospheres in certain (pseudoriemannian) symmetric spaces.
Another is in the Iwasawa decomposition $G = KAN$ of semisimple groups 
with the Killing metric tensor, which need {\em not\/} be (positive or 
negative) definite even on $N$. Here, $K$ is compact and $A$ is abelian.

An early motivation for this study was the observation that there are two
nonisometric pseudoriemannian metrics on the Heisenberg group $H_3$, one
of which is flat.  This is a strong contrast to the Riemannian case in
which there is only one (up to positive homothety) and it is {\em not\/}
flat.  This is not an anomaly, as we now well know.

While the idea of more than one timelike dimension has appeared a few
times in the physics literature, both in string/M-theory and in
brane-world scenarios, essentially all work to date assumes only one.
Thus all applications so far are of Lorentzian or definite nilpotent
groups.  Guediri and co-workers lead the Lorentzian studies, and most of
their results stated near the end of Section \ref{lor} concern a major,
perennial interest in relativity:  the (non)existence of closed timelike
geodesics in compact Lorentzian manifolds.

Others have made use of nilpotent Lie groups with left-invariant (positive
or negative) definite metric tensors, such as Hervig's \cite{H}
constructions of black-hole spacetimes from solvmanifolds (related to
solvable groups:  those with Iwasawa decomposition $G = AN$), including
so-called BTZ constructions.  Definite groups already having received
thorough surveys elsewhere, most notably those of Eberlein, they and their
applications are not included here.

Although the geometric properties of Lie groups with left-invariant
definite metric tensors have been studied extensively, the same has not
occurred for indefinite metric tensors.  For example, while the paper of
Milnor \cite{M} has already become a classic reference, in particular for
the classification of positive definite (Riemannian) metrics on
3-dimensional Lie groups, a classification of the left-invariant
Lorentzian metric tensors on these groups became available only in 1997.
Similarly, only a few partial results in the line of Milnor's study of
definite metrics were previously known for indefinite metrics.  Moreover,
in dimension 3 there are only two types of metric tensors:  Riemannian
(definite) and Lorentzian (indefinite).  But in higher dimensions there
are many distinct types of indefinite metrics while there is still
essentially only one type of definite metric.  This is another reason this
area has special interest now.

As is customary in these volumes, the list of Further Reading consists of 
general survey articles and a select few of the more historically 
important papers. Precise bibliographical information for references 
merely mentioned or alluded to in this article may be found in those. The 
main, general reference on pseudoriemannian geometry is O'Neill's book 
\cite{O}. Eberlein's article \cite{E} covers the Riemannian case. At this 
time, there is no other comprehensive survey of the pseudoriemannian case.
One may use \cite{CP}, \cite {G}, and their reference lists to good 
advantage, however.

\subsection{Inner product and signature}\label{ips}
By an {\em inner product\/} on a vector space $V$ we shall mean a
nondegenerate, symmetric bilinear form on $V$, generally denoted by
$\langle\,,\rangle$.  In particular, we {\em do not\/} assume that it is
positive definite. It has become customary to refer to an ordered pair of 
nonnegative integers $(p,q)$ as the {\em signature\/} of the inner 
product, where $p$ denotes the number of positive eigenvalues and $q$ the 
number of negative eigenvalues. Then nondegeneracy means that $p+q=\dim 
V$.  Note that there is no real geometric difference between $(p,q)$ and 
$(q,p)$; indeed, O'Neill gives handy conversion procedures for this 
(p.\,92f) and for the other major sign variant (curvature, p.\,89). 

A {\em Riemannian\/} inner product has signature $(p,0)$. In view of the 
preceding remark, one might as well regard signature $(0,q)$ as also being
Riemannian, so that ``Riemannian geometry is that of definite metric 
tensors." Similarly, a {\em Lorentzian\/} inner product has either $p=1$
or $q=1$.  In this case, both sign conventions are used in relativistic
theories with the proviso that the ``1" axis is always timelike.

If neither $p$ nor $q$ is 1, there is no physical convention. We shall say
that $v\in V$ is {\em timelike\/} if $\langle v,v\rangle > 0$, {\em
null\/} if $\langle v,v\rangle = 0$, and {\em spacelike\/} if $\langle
v,v\rangle < 0$. (In a Lorentzian example, one may wish to revert to one's 
preferred relativistic convention.) We shall refer to these collectively 
as the {\em causal type\/} of a vector (or of a curve to which a vector is
tangent).

Considering indefinite inner products (and metric tensors) thus 
greatly expands one's purview, from one type of geometry (Riemannian), or 
possibly two (Riemannian and Lorentzian), to a total of $\lfloor
\frac{p+q}{2} \rfloor+1$ distinctly different types of geometries on {\em
the same underlying differential manifolds.}

\subsection{Rise of 2-step groups}\label{2step}
Throughout, $N$ will denote a connected (and simply connected, usually),
nilpotent Lie group with Lie algebra $\n$ having center $\z$.  We shall
use $\langle\,,\rangle$ to denote either an inner product on $\n$ or the
induced left-invariant pseudoriemannian (indefinite) metric tensor on $N$.

For all nilpotent Lie groups, the exponential map $\exp :\n\to N$
is surjective. Indeed, it is a diffeomorphism for simply connected $N$; in
this case we shall denote the inverse by $\log$.

One of the earliest papers on the Riemannian geometry of nilpotent Lie 
groups is \cite{W}. Since then, a few other papers about general nilpotent
Lie groups have appeared, including \cite{K} and \cite{P}, but the area 
has not seen a lot of progress.

In 1981, Kaplan published \cite{K1} and everything changed. Following this 
paper and its successor \cite{K2}, almost all subsequent work on the
left-invariant geometry of nilpotent groups has been on 2-step groups.

Briefly, Kaplan defined a new class of nilpotent Lie groups, calling them 
{\em of Heisenberg type.} This was soon abbreviated to $H$-type, and has 
since been called Heisenberg-like and (unfortunately) generalized 
Heisenberg. (Unfortunate, because that term was already in use for another
class, not all of which are of $H$-type.) What made them so compelling was
that (almost) everything was {\em explicitly\/} calculable, thus making 
them the next great test bed after symmetric spaces. 

\begin{definition}
We say that $N$ (or $\n$) is 2-step nilpotent when $[\n,\n]\subseteq\z$.
\end{definition}
Then $[[\n,\n],\n] = 0$ and the generalization to $k$-step nilpotent is 
clear:
$$[[\cdots[[[\n,\n],\n],\n]\cdots],\n]=0$$
with $k+1$ copies of $\n$ (or $k$ nested brackets, if you prefer).

It soon became apparent that $H$-type groups comprised a subclass of
2-step groups; for a nice, modern proof see \cite{BTV}.  By around 1990,
they had also attracted the attention of the spectral geometry industry,
and Eberlein produced the seminal survey (with important new results) from
which the modern era began.  (It was published in 1994 \cite{E'}, but the
preprint had circulated widely since 1990.) Since then, activity around 
2-step nilpotent Lie groups has mushroomed; see the references in \cite{E}.

Finally, turning to pseudoriemannian nilpotent Lie groups, with perhaps
one or two exceptions, all results so far have been obtained only for
2-step groups. Thus the remaining sections of this article will be devoted
almost exclusively to them.

The Baker-Campbell-Hausdorff formula takes on a particularly simple form 
in these groups:
\begin{equation}
\exp(x)\exp(y)=\exp(x+y+\half[x,y])\, .
\end{equation}
\begin{proposition}
In a pseudoriemannian 2-step nilpotent Lie group, the exponential map
preserves causal character.  Alternatively, 1-parameter subgroups are
curves of constant causal character.
\end{proposition}
Of course, 1-parameter subgroups need not be geodesics.

\subsection{Lattices and completeness}\label{lac}
We shall need some basic facts about lattices in $N$.  In nilpotent
Lie groups, a lattice is a discrete subgroup $\G$ such that the
homogeneous space $M=\G\bs N$ is compact. Here we follow the convention 
that a lattice acts on the {\em left,} so that the coset space consists of
{\em left\/} cosets and this is indicated by the notation. Other subgroups
will generally act on the {\em right,} allowing better separation of the 
effects of two simultaneous actions.

Lattices do not always exist in nilpotent Lie groups.
\begin{theorem}
The simply connected, nilpotent Lie group $N$ admits a lattice if and only 
if there exists a basis of its Lie algebra $\n$ for which the structure 
constants are rational.
\end{theorem}
Such a group is said to have a rational structure, or simply to
be rational. 

A {\em nilmanifold\/} is a (compact) homogeneous space of the form $\G\bs 
N$ where $N$ is a connected, simply connected (rational) nilpotent Lie
group and $\G$ is a lattice in $N$. An {\em infranilmanifold\/} has a
nilmanifold as a finite covering space. They are commonly regarded as a 
noncommutative generalization of tori, the Klein bottle being the simplest
example of an infranilmanifold that is not a nilmanifold.

We recall the result of Marsden from \cite{O}.
\begin{theorem}
A compact, homogeneous pseudoriemannian space is geodesically complete.
\end{theorem}
Thus if a rational $N$ is provided with a bi-invariant metric tensor 
$\langle\, ,\rangle$, then $M$ becomes a compact, homogeneous 
pseudoriemannian space which is therefore complete. It follows that 
$(N,\langle\, ,\rangle)$ is itself complete. In general, however, the
metric tensor is not bi-invariant and $N$ need not be complete.

For 2-step nilpotent Lie groups, things work nicely as shown by this
result first published by Guediri.
\begin{theorem}
On a 2-step nilpotent Lie group, all left-invariant pseudoriemannian 
metrics are geodesically complete.
\end{theorem}
No such general result holds for 3- and higher-step groups, however.

\section{2-step Groups}\label{2sg}
In the Riemannian (positive-definite) case, one splits $\n = \z\ds\v =
\z\ds\z^{\perp}$ where the superscript denotes the orthogonal complement
with respect to the inner product $\langle\,,\rangle$.  In the general
pseudoriemannian case, however, $\z\ds\z^{\perp}\not=\n$.  The problem is
that $\z$ might be a {\em degenerate\/} subspace; {\em i.e.,} it might
contain a null subspace $\U$ for which $\U\subseteq\U^{\perp}$.

It turns out that this possible degeneracy of the center causes the
essential differences between the Riemannian and pseudoriemannian cases.
So far, the only general success in studying groups with degenerate
centers was in \cite{CP} where we used an adapted Witt decomposition of
$\n$ together with an involution $\iota$ exchanging the two null parts.

Observe that if $\z$ is degenerate, the null subspace $\U$ is well defined
invariantly. We shall use a decomposition
\begin{equation}
\n = \z\ds\v = \U\ds\Z\ds\V\ds\E
\end{equation}
in which $\z = \U\ds\Z$ and $\v = \V\ds\E$, $\U$ and $\V$ are complementary
null subspaces, and $\U^{\perp}\cap\V^{\perp} = \Z\ds\E$. Although the
choice of $\V$ is {\em not\/} well defined invariantly, once a $\V$ has been 
chosen then $\Z$ and $\E$ {\em are\/} well defined invariantly. Indeed, $\Z$
is the portion of the center $\z$ in $\U^{\perp}\cap\V^{\perp}$ and $\E$
is its orthocomplement in $\U^{\perp}\cap\V^{\perp}$.  This is a Witt
decomposition of $\n$ given $\U$, easily seen by noting that $\left(\U \ds
\V\right)^{\perp} = \Z\ds\E$, adapted to the special role of the center in
$\n$.

We shall also need to use an involution $\io$ that interchanges $\U$ and
$\V$ and which reduces to the identity on $\Z\oplus\E$ in the Riemannian
(positive-definite) case.  (The particular choice of such an involution is
not significant.)  It turns out that $\io$ is an isometry of $\n$ which
does {\em not\/} integrate to an isometry of $N$.  The adjoint with
respect to $\langle\,,\rangle$ of the adjoint representation of the Lie
algebra $\n$ on itself is denoted by $\add{}{}\!$.
\begin{definition}
The linear mapping
$$ j:\U\ds\Z\rightarrow\End\left(\V\ds\E\right) $$
is given by 
$$ j(a)x = \io\add{x}{\io a}\, .$$
\end{definition}

Formulas for the connection and curvatures, and explicit forms for many
examples, may be found in \cite{CP}.  It turns out there is a relatively
large class of flat spaces, a clear distinction from the Riemannian case
in which there are none.  

Let $x,y\in\n$.  Recall that {\em homaloidal\/} planes are those for which
the numerator $\langle R(x,y)y,x\rangle$ of the sectional curvature
formula vanishes.  This notion is useful for degenerate planes tangent to
spaces that are not of constant curvature.
\begin{definition}
A submanifold of a pseudoriemannian manifold is {\em flat\/}
if and only if every plane tangent to the submanifold is homaloidal.
\end{definition}
\begin{theorem}
The center $Z$ of $N$ is flat.
\end{theorem}
\begin{corollary}
The only $N$ of constant curvature are flat.
\end{corollary}

The degenerate part of the center can have a profound effect on the
geometry of the whole group.
\begin{theorem}
If\/ $[\n,\n] \subseteq \U$ and\label{e0f} $\E=\{0\}$, then $N$ is flat.
\end{theorem}
Among these spaces, those that also have  $\Z = \{0\}$ (which condition
itself implies $[\n,\n] \subseteq \U$) are fundamental, with the more
general ones obtained by making nondegenerate central extensions.  It is
also easy to see that the product of any flat group with a nondegenerate
abelian factor is still flat.

This is the best possible result in general. Using weaker hypotheses in
place of $\E = \{0\}$, such as $[\V,\V] = \{0\} = [\E,\E]$, it is easy to
construct examples which are not flat.
\begin{corollary}
If $\dim Z \ge \lceil\frac{n}{2}\rceil$, then there exists a flat metric
on $N$.
\end{corollary}
Here $\lceil r\rceil$ denotes the least integer greater than or equal to
$r$ and $n=\dim N$.

Before continuing, we pause to collect some facts about the condition
$[\n,\n] \subseteq \U$ and its consequences.
\begin{remark}
Since it implies $j(z) = 0$ for all $z\in\Z$, this latter is possible with
no pseudoeuclidean de Rham factor, unlike the Riemannian case.  (On the
other hand, a pseudoeuclidean de Rham factor is characterized in terms of
the Kaplan-Eberlein map $j$ whenever the center is nondegenerate.)

Also, it implies $j(u)$ interchanges $\V$ and $\E$ for all $u\in\U$ if and
only if $[\V,\V] = [\E,\E] = \{0\}$.  Examples are the Heisenberg group
and the groups $H(p,1)$ for $p\ge 2$ with null centers.

Finally we note it implies that, for every $u\in\U$, $j(u)$ maps $\V$ to 
$\V$ if and only if $j(u)$ maps $\E$ to $\E$ if and only if $[\V,\E] = \{0\}$.
\end{remark}
\begin{proposition}
If $j(z) = 0$ for all $z\in\Z$ and $j(u)$ interchanges $\V$
and $\E$ for all $u\in\U$, then $N$ is Ricci flat.
\end{proposition}
\begin{proposition}
If $j(z) = 0$ for all $z\in\Z$, then $N$ is scalar flat.  In particular,
this occurs when $[\n,\n]\subseteq\U$.
\end{proposition}

Much like the Riemannian case, we would expect that $(N,\langle\,
,\rangle)$ should in some sense be similar to flat pseudoeuclidean space.
This is seen, for example, {\em via\/} the existence of totally geodesic
subgroups \cite{CP}. (O'Neill \cite[Ex.\,9, p.\,125]{O} has extended the
definition of totally geodesic to degenerate submanifolds of
pseudoriemannian manifolds.)
\begin{example}
For any $x\in\n$ the 1-parameter subgroup $\exp (tx)$ is a geodesic if and
only if $x\in\z$ or $x\in\U\ds\E$.  This is essentially the same as the
Riemannian case, but with some additional geodesic 1-parameter subgroups
coming from $\U$.
\end{example}
\begin{example}
Abelian subspaces of $\V\ds\E$ are Lie subalgebras of $\n$, and give rise
to complete, flat, totally geodesic abelian subgroups of $N$, just as in
the Riemannian case.  Eberlein's construction is valid in general, and
shows that if $\dim\V \ds\E \geq 1 + k + k\dim\z$, then every nonzero
element of $\V\ds\E$ lies in an abelian subspace of dimension $k+1$.
\end{example}
\begin{example}
The center $Z$ of $N$ is a complete, flat, totally geodesic submanifold.
Moreover, it determines a foliation of $N$ by its left translates, so each
leaf is flat and totally geodesic, as in the Riemannian case.  In the
pseudoriemannian case, this foliation in turn is the orthogonal direct sum
of two foliations determined by $\U$ and $\Z$, and the leaves of the
$\U$-foliation are also null.  All these leaves are complete.
\end{example}

There is also the existence of $\dim\z$ independent
first integrals, a familiar result in pseudoeuclidean space, and the
geodesic equations are completely integrable; in certain cases (mostly
when the center is nondegenerate) one can obtain {\em explicit\/}
formulas.  Unlike the Riemannian case, there are flat groups 
(nonabelian) which are isometric to pseudoeuclidean spaces (abelian).
\begin{theorem}
If $[\n,\n] \subseteq \U$ and $\E=\{0\}$, then $N$ is geodesically
connected.  Consequently, so is any nilmanifold with such a universal
covering space.
\end{theorem}
Thus these compact nilmanifolds are much like tori.  This is also
illustrated by the computation of their period spectrum.

\subsection{Isometry group}\label{isog}
The main new feature is that when the center is degenerate, the isometry
group can be strictly larger in a significant way than when the center is
nondegenerate (which includes the Riemannian case).

Letting $\Aut(N)$ denote the automorphism group of $N$ and $I(N)$ the
isometry group of $N$, set $O(N) = \Aut(N)\cap I(N)$.  In the Riemannian
case, $I(N) = O(N)\ltimes N$, the semidirect product where $N$ acts as
left translations.  We have chosen the notation $O(N)$ to suggest an
analogy with the pseudoeuclidean case in which this subgroup is precisely
the (general, including reflections) pseudorthogonal group.  According to
Wilson \cite{Wi}, this analogy is good for {\em any\/} nilmanifold (not
necessarily 2-step).

To see what is true about the isometry group in general, first consider the
(left-invariant) splitting of the tangent bundle $TN = \z N \ds \v N$.
\begin{definition}
Denote by $\Ispl(N)$ the subgroup of the isometry group $I(N)$ which
preserves the splitting $TN = \z N \ds \v N$.  Further, let $\Iaut(N) =
O(N)\ltimes N$, where $N$ acts by left translations.
\end{definition}
\begin{proposition}
If $N$ is a simply-connected, 2-step nilpotent Lie group with
left-invariant metric tensor, then $\Ispl(N) \le \Iaut(N)$.
\end{proposition}
There are examples to show that $\Ispl < \Iaut$ is possible when $\U 
\ne \{0\}$.

When the center is degenerate, the relevant group analogous to a
pseudorthogonal group may be larger.
\begin{proposition}
Let $\wO(N)$ denote the subgroup of $I(N)$ which fixes $1 \in N$.  Then
$I(N) \cong \wO(N) \ltimes N$, where $N$ acts by left translations.
\end{proposition}
The proof is obvious from the definition of $\wO$.  It is also obvious that 
$O \le \wO$.  Examples show that $O < \wO$, hence $\Iaut < I$, is possible
when the center is degenerate.

Thus we have three groups of isometries, not necessarily equal in general:
$\Ispl \le \Iaut \le I$.  When the center is nondegenerate ($\U = \{0\}$),
the Ricci transformation is block-diagonalizable and the rest of Kaplan's
proof using it now also works.
\begin{corollary}
If the center is nondegenerate, then $I(N) = \Ispl(N)$, hence $\wO(N)
\cong O(N)$.
\end{corollary}

In the next few results, we use the phrase ``a subgroup isometric to" a 
group to mean that the isometry is also an isomorphism of groups.
\begin{proposition}
For any $N$ containing a subgroup isometric to the flat 3-dimensional 
Heisenberg group, $$\Ispl(N) < \Iaut(N) < I(N) .$$
\end{proposition}
Unfortunately, this class does not include our flat groups in which
$[\n,\n] \subseteq \U$ and $\E = \{0\}$.  However, it does include many
groups that do not satisfy $[\n,\n] \subseteq \U$, such as the simplest
quaternionic Heisenberg group.
\begin{remark}
A direct computation shows that on this flat $H_3$ with null center, the
only Killing fields with geodesic integral curves are the nonzero scalar
multiples of a vector field tangent to the center.
\end{remark}
\begin{proposition}
For any $N$ containing a subgroup isometric to the flat $H_3\times\R$ 
with null center,
$$\Ispl(N) < \Iaut(N) < I(N) .$$
\end{proposition}
Many of our flat groups in which $[\n,\n] \subseteq \U$ and $\E = \{0\}$
have such a subgroup isometrically embedded, as in fact do many others
which are not flat.

\subsection{Lattices and periodic geodesics}\label{lpg}
In this subsection, we assume that $N$ is rational and let $\G$
be a lattice in $N$.

Certain tori $T_F$ and $T_B$ provide the model fiber and the base for a
submersion of the coset space $\G\bs N$.  This submersion may not be
pseudoriemannian in the usual sense, because the tori may be degenerate.
We began the study of periodic geodesics in these compact nilmanifolds,
and obtained a complete calculation of the period spectrum for certain 
flat spaces.

To the compact nilmanifold $\G\bs N$ we may associate two flat 
(possibly degenerate) tori.
\begin{definition}
Let $N$ be a simply connected, 2-step nilpotent Lie group with lattice 
$\G$ and let $\pi:\n\to\v$ denote the projection. Define
\begin{eqnarray*}
T_{\z} &=& \z/(\log\G\cap\z)\, ,\\
T_{\v} &=& \v/\pi(\log\G)\, .
\end{eqnarray*}
\end{definition}
Observe that $\dim T_{\z} + \dim T_{\v} = \dim\z + \dim\v = \dim\n$.

Let $m=\dim\z$ and $n=\dim\v$.  It is a consequence of a theorem of Palais
and Stewart that $\G\bs N$ is a principal $T^m$-bundle over $T^n$.  The
model fiber $T^m$ can be given a geometric structure from its closed
embedding in $\G\bs N$; we denote this geometric $m$-torus by $T_F$.
Similarly, we wish to provide the base $n$-torus with a geometric
structure so that the projection $p_B:\G\bs N\surj T_B$ is the appropriate
generalization of a pseudoriemannian submersion \cite{O} to (possibly)
degenerate spaces.  Observe that the splitting $\n = \z\ds\v$ induces
splittings $TN = \z N\ds\v N$ and $T(\G\bs N) = \z(\G\bs N)\ds\v(\G\bs
N)$, and that $p_{B*}$ just mods out $\z(\G\bs N)$.  Examining O'Neill's
definition, we see that the key is to construct the geometry of $T_B$ by
defining
\begin{equation}
p_{B*} : \v_\eta(\G\bs N)\to T_{p_B(\eta)}(T_B) \mbox{ for each }
\eta\in\G\bs N \mbox{ is an isometry }
\end{equation}
and
\begin{equation}
\nabla^{T_B}_{p_{B*}x}p_{B*}y = p_{B*}\left(\pi\nabla_x y\right)
\mbox{ for all $x,y\in\v=\V\ds\E$, }
\end{equation}
where $\pi:\n\to\v$ is the projection.  Then the rest of the usual results
will continue to hold, provided that sectional curvature is replaced by
the numerator of the sectional curvature formula at least when elements of
$\V$ are involved:
\begin{equation}
\langle R_{T_B}(p_{B*}x,p_{B*}y)p_{B*}y,p_{B*}x\rangle = \langle
R_{\Gamma\bs N}(x,y)y,x\rangle + {\txs\frac{3}{4}}
\langle[x,y],[x,y]\rangle .
\end{equation}
Now $p_B$ will be a pseudoriemannian submersion in the usual sense if and
only if $\U=\V=\{0\}$, as is always the case for Riemannian spaces.

In the Riemannian case, Eberlein showed that $T_F\cong T_{\z}$ and
$T_B\cong T_{\v}$.  In general, $T_B$ is flat only if $N$ has a
nondegenerate center or is flat.
\begin{remark}
Observe that the torus $T_B$ may be decomposed into a topological product
$T_E \times T_V$ in the obvious way.  It is easy to check that $T_E$ is
flat and isometric to $(\log\G \cap \E)\bs\E$, and that $T_V$ has a linear
connection not coming from a metric and not flat in general.  Moreover,
the geometry of the product is ``twisted" in a certain way.  It would be
interesting to determine which tori could appear as such a $T_V$ and how.
\end{remark}
\begin{theorem}
Let $N$ be a simply connected, 2-step nilpotent Lie group with lattice
$\G$, a left-invariant metric tensor, and tori as above.  The fibers $T_F$
of the (generalized) pseudoriemannian submersion $\G\bs N\surj T_B$ are
isometric to $T_{\z}$.  If in addition the center $Z$ of $N$ is
nondegenerate, then the base $T_{B}$ is isometric to $T_{\v}$.
\end{theorem}

We recall that elements of $N$ can be identified with elements of the
isometry group $I(N)$: namely, $n\in N$ is identified with the isometry
$\phi = L_n$ of left translation by $n$. We shall abbreviate this by
writing $\phi\in N$.
\begin{definition}
We say that $\phi\in N$ {\em translates\/} the geodesic $\g$ by $\om$ if
and only if $\phi\g(t) = \g(t+\om)$ for all $t$.  If $\g$ is a unit-speed
geodesic, we say that $\om$ is a {\em period\/} of $\phi$.
\end{definition} 
Recall that unit speed means that $|\dot{\g}| = \left|\langle\dot{\g},
\dot{\g}\rangle\right|^{\frac{1}{2}} = 1$.  Since there is no natural
normalization for null geodesics, we do not define periods for them.  In
the Riemannian case and in the timelike Lorentzian case in strongly causal
spacetimes, unit-speed geodesics are parameterized by arclength and this
period is a translation distance.  If $\phi$ belongs to a lattice $\G$, it
is the length of a closed geodesic in $\G\bs N$.

In general, recall that if $\g$ is a geodesic in $N$ and if $p_N:N\surj
\G\bs N$ denotes the natural projection, then $p_N\g$ is a periodic
geodesic in $\G\bs N$ if and only if some $\phi\in\G$ translates $\g$.  We
say {\em periodic\/} rather than {\em closed\/} here because in
pseudoriemannian spaces it is possible for a null geodesic to be closed
but not periodic.  If the space is geodesically complete or Riemannian,
however, then this does not occur; the former is in fact the case for our
2-step nilpotent Lie groups.  Further recall that free homotopy classes of
closed curves in $\G\bs N$ correspond bijectively with conjugacy classes
in $\G$.
\begin{definition}
Let $\sC$ denote either a nontrivial, free homotopy class of
closed curves in $\G\bs N$ or the corresponding conjugacy class in
$\G$. We define $\wp(\sC)$ to be the set of all periods of periodic
unit-speed geodesics that belong to $\sC$.
\end{definition}
In the Riemannian case, this is the set of lengths of closed geodesics in
$\sC$, frequently denoted by $\ell(\sC)$.
\begin{definition}
The {\em period spectrum\/} of $\G\bs N$ is the set 
$$ \specp(\G\bs N) = \bigcup_{\sC}\wp(\sC)\,,$$
where the union is taken over all nontrivial, free homotopy classes of
closed curves in $\G\bs N$.
\end{definition}
In the Riemannian case, this is the length spectrum $\specl(\G\bs N)$.
\begin{example}
Similar to the Riemannian case, we can compute the period spectrum of a
flat torus $\G\bs\R^m$, where $\G$ is a lattice (of maximal rank,
isomorphic to ${\Bbb Z}^m$).  Using calculations in an analogous way as
for finding the length spectrum of a Riemannian flat torus, we easily
obtain
$$ \specp(\G\bs\R^m) = \{ |g|\ne 0 \mid g\in\G\}\,. $$
It is also easy to see that the nonzero d'Alembertian spectrum is related
to the analogous set produced from the dual lattice $\G^*$ as multiples by
$\pm 4\pi^2$, almost as in the Riemannian case.
\end{example}

As in this example, simple determinacy of periods of unit-speed geodesics
helps make calculation of the period spectrum possible purely in terms of
$\log\G \subseteq \n$.

For the rest of this subsection, we assume that $N$ is a simply connected,
2-step nilpotent Lie group with left-invariant pseudoriemannian metric
tensor $\langle\,,\rangle$.  Note that non-null geodesics may be taken to
be of unit speed.  Most nonidentity elements of $N$ translate some
geodesic, but not necessarily one of unit speed.

For our special class of flat 2-step nilmanifolds, we can calculate the
period spectrum completely.
\begin{theorem}
If $[\n,\n] \subseteq \U$ and $\E=\{0\}$, then $\specp(M)$ can be
completely calculated from $\log\G$ for any $M = \G\bs N$.
\end{theorem}
Thus we see again just how much these flat, 2-step nilmanifolds are like
tori.  All periods can be calculated purely from $\log\G \subseteq \n$,
although some will not show up from the tori in the fibration.
\begin{corollary}
$\specp(T_B)$ (respectively, $T_F$) is $\cup_{\sC}\,\wp(\sC)$ where the
union is taken over all those free homotopy classes\/ $\sC$ of closed
curves in $M = \G\bs N$ that\/ {\em do not} (respectively,\/ {\em do})
contain an element in the center of\/ $\G \cong \pi_1(M)$, except for
those periods arising only from unit-speed geodesics in $M$ that project
to null geodesics in both $T_B$ and $T_F$.
\end{corollary}
We note that one might consider using this to assign periods to some null
geodesics in the tori $T_B$ and $T_F$.

When the center is nondegenerate, we obtain results similar to Eberlein's.
Here is part of them.
\begin{theorem}
Assume $\U = \{0\}$.  Let $\phi \in N$ and write $\log\phi = z^* + e^*$.
Assume $\phi$ translates the unit-speed geodesic $\g$ by $\om > 0$.  Let
$z'$ denote the component of $z^*$ orthogonal to $[e^*,\n]$ and set $\om^*
= |z' + e^*|$. Let $\dg(0) = z_0 + e_0$. Then
\begin{enumerate}
\item $|e^*| \le \om$. In addition, $\om < \om^*$ for timelike
   (spacelike) geodesics with $\om z_0 - z'$ timelike (spacelike), and
   $\om > \om^*$ for timelike (spacelike) geodesics with $\om z_0 - z'$
   spacelike (timelike).

\item $\om = |e^*|$ if and only if $\g(t) = \exp\left(t\,e^*\!/|e^*|
   \right)$ for all $t \in \R$.

\item $\om = \om^*$ if and only if $\om z_0 - z'$ is null.
\end{enumerate}
\end{theorem}
Although $\om^*$ need not be an upper bound for periods as in the
Riemannian case, it nonetheless plays a special role among all periods, as
seen in item 3 above, and we shall refer to it as the {\em
distinguished\/} period associated with $\phi \in N$.  When the center is
definite, for example, we do have $\om \le \om^*$.

Now the following definitions make sense at least for $N$ with a
nondegenerate center.
\begin{definition}
Let $\sC$ denote either a nontrivial, free homotopy class of closed curves
in $\G\bs N$ or the corresponding conjugacy class in $\G$.  We define
$\wp^*(\sC)$ to be the distinguished periods of periodic unit-speed
geodesics that belong to $\sC$.
\end{definition}
\begin{definition}
The {\em distinguished period spectrum\/} of $\G\bs N$ is the set
$$ \Dspecp(\G\bs N) = \bigcup_{\sC}\wp^*(\sC)\,,$$
where the union is taken over all nontrivial, free homotopy classes of
closed curves in $\G\bs N$.
\end{definition}
Then we get this result.
\begin{corollary}
Assume the center is nondegenerate.  If $\n$ is nonsingular, then
$\specp(T_B)$ (respectively, $T_F$) is precisely the period spectrum
(respectively, the distinguished period spectrum) of those free homotopy
classes\/ $\sC$ of closed curves in $M = \G\bs N$ that\/ {\em do not}
(respectively,\/ {\em do}) contain an element in the center of\/ $\G \cong
\pi_1(M)$, except for those periods arising only from unit-speed geodesics
in $M$ that project to null geodesics in both $T_B$ and $T_F$.
\end{corollary}

\subsection{Conjugate loci}\label{cjl}
This is the only general result on conjugate points.
\begin{proposition}
Let $N$ be a simply connected, 2-step nilpotent Lie group with
left-invariant metric tensor $\langle\,,\rangle$, and let $\g$ be a
geodesic with $\dg(0) = a\in \z$.  If\/ $\add{\bul}{a} = 0$, then there
are no conjugate points along $\g$.
\end{proposition}
In the rest of this subsection, we assume that the center of $N$ is 
nondegenerate.

For convenience, we shall use the notation $J_z = \add{\bul}{z}$ for any
$z\in\z$. (Since the center is nondegenerate, the involution $\io$ may
be omitted.) We follow Ciatti \cite{C} for this next definition. As in the 
Riemannian case, one might as well make 2-step nilpotent part of the 
definition since it effectively is so anyway.
\begin{definition}
$N$ is said to be of {\em pseudo$H$-type} if and only if
$$ J^2_z = -\langle z, z \rangle I $$
for any $z \in \z$.
\end{definition}
Complete results on conjugate loci have been obtained only for these 
groups \cite{JP}. For example, using standard results from analytic 
function theory, one can show that the conjugate locus is an analytic 
variety in $N$. This is probably true for general 2-step groups, but the 
proof I know works only for pseudo$H$-type.

\begin{definition}
Let $\g$ denote a geodesic and assume that $\g(t_0)$ is conjugate to 
$\g(0)$ along $\g$. To indicate that the multiplicity of $\g(t_0)$ is $m$,
we shall write $\mcp(t_0)=m$. To distinguish the notions clearly, we shall
denote the multiplicity of $\l$ as an eigenvalue of a specified linear 
transformation by $\mev{\l}$.
\end{definition}
Let $\g$ be a geodesic with $\g(0)=1$ and $\dg(0) = z_0+x_0 \in
\z\ds\v$, respectively, and let $J = J_{z_0}$. If $\g$ is not
null, we may assume $\g$ is normalized so that $\langle\dg,\dg\rangle =
\pm 1$. As usual, $\bZ^*$ denotes the set of all integers with 0 removed.
\begin{theorem}
Under these assumptions, if $N$ is of pseudo$H$-type then:
\begin{enumerate}
\item if $z_0=0$ and $x_0\neq 0$, then $\g(t)$ is conjugate to $\g(0)$
along $\g$ if and only if $\langle x_0,x_0\rangle < 0$ and
$$-\frac{12}{t^2} = \langle x_0,x_0\rangle\, ,$$
in which case $\mcp(t) = \dim\z$;
\item if $z_0\neq 0$ and $x_0 = 0$, then $\g(t)$ is conjugate to $\g(0)$
along $\g$ if and only if $\langle z_0,z_0\rangle > 0$ and
$$t \in  \frac{2\pi}{|z_0|} \bZ^* ,$$
in which case $\mcp(t) = \dim\v$.
\end{enumerate}
\end{theorem}
\begin{theorem}
Let $\g$ be such a geodesic in a pseudo$H$-type group $N$ with
$z_0 \neq 0 \neq x_0$.
\begin{enumerate}
\item If $\langle z_0 , z_0 \rangle =\alpha^2$ with $\alpha>0$, then
$\gamma (t_0)$ is conjugate to $\g(0)$ along $\gamma$ if and only if
$$ t_0 \in \frac{2\pi}{\alpha }\bZ^* \cup \A_1 \cup \A_2$$
where
$$ \A_1 = \left\{ t \in \R \biggm| \langle x_0 , x_0 \rangle\frac{\alpha
   t}{2}\cot\frac{\alpha t}{2} = \langle \dg , \dg \rangle \right\}
   \quad\mbox{and \qquad } $$
$$ \A_2 = \left\{ t\in \R \biggm| \alpha t = \frac{\langle x_0 , x_0
   \rangle }{\langle \dg , \dg \rangle + \langle z_0 , z_0 \rangle }
   \sin\alpha t\right\}\quad\mbox{when }\,\dim\z  \geq 2\,. $$
If $t_0\in (2\pi/\alpha)\bZ^*$, then
$$ \mcp(t_0) = \left\{ \begin{array}{cl}
\dim\v-1 & \mbox{ if }\,\langle\dg,\dg\rangle + \langle z_0,z_0\rangle \ne
   0\,,\\
\dim\n-2 & \mbox{ if }\,\langle\dg,\dg\rangle + \langle z_0,z_0\rangle = 
   0\,. \end{array} \right. $$
If $t_0\notin (2\pi/\alpha)\bZ^*$, then
$$ \mcp(t_0) = \left\{ \begin{array}{cl}
1 & \mbox{ if }\, t_0\in\A_1 - \A_2\,,\\
\dim\z-1 & \mbox{ if }\, t_0\in\A_2 - \A_1\,,\\
\dim\z & \mbox{ if }\, t_0\in\A_1\cap\A_2\,. \end{array} \right. $$

\item If $\langle z_0 , z_0 \rangle = -\beta^2$ with $\beta >0$ , then
$\gamma (t_0)$ is a conjugate point along $\gamma$ if and only if $t_0 \in
\B_1 \cup \B_2$ where
$$ \B_1 = \left\{ t \in \R \biggm| \langle x_0 , x_0 \rangle \frac{\beta
   t}{2}\coth\frac{\beta t}{2} = \langle \dg , \dg \rangle \right\}
   \quad\mbox{and \qquad } $$
$$ \B_2 = \left\{ t\in \R \biggm| \beta t = \frac{\langle x_0 , x_0
   \rangle }{\langle \dg , \dg \rangle + \langle z_0 , z_0 \rangle } \sinh
   \beta t\right\}\quad\mbox{when }\,\dim\z \geq 2\,. $$
The multiplicity is
$$ \mcp(t_0) = \left\{ \begin{array}{cl}
1 & \mbox{ if }\, t_0\in\B_1 - \B_2\,,\\
\dim\z-1 & \mbox{ if }\, t_0\in\B_2 - \B_1\,,\\
\dim\z & \mbox{ if }\, t_0\in\B_1\cap\B_2\,. \end{array} \right. $$

\item If $ \langle z_0 , z_0 \rangle = 0$, then $ \gamma (t_0) $ is a
conjugate point along $\gamma$ if and only if
$$ t_0^2 = -\frac{12}{\langle x_0 , x_0 \rangle }\, ,$$
and $\mcp(t_0)=\dim\z-1$.

\end{enumerate}
\end{theorem}
This covers all cases for a pseudo$H$-type group with a center of any
dimension.

Some results on other 2-step groups and examples (including pictures in 
dimension 3) may be found in the references cited in \cite{JP}. When the 
groups are not pseudo$H$-type, however, complete results are available 
only when the center is 1-dimensional. Guediri \cite{G'} has results in 
the timelike Lorentzian case.

\section{Lorentzian Groups}\label{lor}
Not too long ago, only a few partial results in the line of Milnor's study
of definite metrics were known for indefinite metrics \cite{B,N}, and 
they were Lorentzian.

Guediri \cite{G} and others have made special study of Lorentzian 2-step 
groups, partly because of their relevance to General Relativity where they
can be used to provide interesting and important (counter)examples. 
Special features of Lorentzian geometry frequently enable them to obtain
much more complete and explicit results than are possible in general.

For example, in \cite{G} Guediri was able to provide a complete and 
explicit integration of the geodesic equations for Lorentzian 2-step 
groups. This includes the case of a degenerate center, which only required 
extremely careful handling through a number of cases. He also paid special
attention to the existence of closed timelike geodesics, reflecting the 
relativistic concerns.

As usual, $N$ denotes a connected and simply connected 2-step nilpotent 
Lie group. For the rest of this section we assume the left-invariant metric
tensor is Lorentzian. Whenever a lattice is mentioned, we also assume the 
group is rational.
\begin{proposition}
If the center is degenerate, then no timelike geodesic can be translated 
by a central element.
\end{proposition}
Thus, there can be no closed timelike geodesics parallel to the center in 
any nilmanifold obtained from such an $N$.
\begin{theorem}
If the center is Lorentzian, then $\G\bs N$ contains no timelike or null 
closed geodesics for any lattice $\G$.
\end{theorem}

To handle degenerate centers, three refined notions for nonsingular are
used:  {\em almost,} {\em weakly,} and {\em strongly\/} nonsingular.
The precise definitions involve an adapted Witt decomposition (as in the
general pseudoriemannian case, but a rather different one here) and are
quite technical, as is typical.  We refer to \cite{G} for details.
\begin{theorem}
If $N$ is weakly nonsingular, then no timelike geodesic can be translated 
by an element of $N$.
\end{theorem}
\begin{corollary}
If $N$ is flat, then no timelike geodesic can be translated by a 
nonidentity element.
\end{corollary}
\begin{corollary}
If $N$ is flat, then $\G\bs N$ contains no closed timelike geodesics for 
any lattice $\G$.
\end{corollary}
\begin{corollary}
If $N$ is weakly nonsingular, then $\G\bs N$ contains no closed timelike 
geodesic.
\end{corollary}
\begin{corollary}
If $H_{2k+1}$ is a Lorentzian Heisenberg group with degenerate center, 
then $\G\bs N$ contains no closed timelike geodesic.
\end{corollary}

Guediri also has the only non-Riemannian results so far about the 
phenomenon Eberlein called ``in resonance." Roughly speaking, this occurs 
when the eigenvalues of the map $j$ have rational ratios. (The Lorentzian 
case actually requires a slightly more complicated condition when the
center is degenerate.)
\begin{theorem}
If $N$ is almost nonsingular, then $N$ is in resonance if and only if 
every geodesic of $N$ is translated by some element of $N$.
\end{theorem}

\end{document}